\documentclass [12pt]{amsart}

\usepackage[left=1in,right=1in,top=1in,bottom=1in]{geometry} 

\usepackage{amsmath}

\usepackage{url}

\usepackage{amssymb,amsthm,amsmath}

\usepackage{ifthen}

\usepackage{graphicx}

\DeclareGraphicsExtensions{.eps}

\newtheorem{theorem}{Theorem}[section]

\theoremstyle{definition}

\newtheorem{definition}[theorem]{Definition}

\newtheorem{example}[theorem]{Example}

\newtheorem{remark}[theorem]{Remark}

\numberwithin{equation}{section}

\begin{document}
	
	\date{\today}
	\title{Operator version of Korovkin theorem; Degree of convergence and its applications
 }
	\author{V. B. Kiran Kumar}
	\email{vbk@cusat.ac.in} 
	\author{P. C. Vinaya}
	\email{vinayapc01@gmail.com } 
	
	\address{Department of Mathematics, Cochin University of Science And Technology, Kerala, India.}
	\begin{abstract}
		In a recent article, Dumitru Popa proved an operator version of the Korovkin theorem. We recall the quantitative version of the Korovkin theorem obtained by O. Shisha and B. Mond in $1968$. In this paper, we obtain a quantitative estimate for the operator version of the Korovkin theorem obtained by Dumitru Popa. We also consider various examples where the operator version is applicable and obtain similar estimates leading to the degree of convergence. In addition, we obtain the trigonometric analogue of this result by proving the quantitative version. Finally, we apply this result to the preconditioning problem of large linear systems with the Toeplitz structure.
		\end{abstract}
	\maketitle	
\section{Introduction}
	The discovery of the famous Korovkin's theorem in $1953$ by P.P. Korovkin was a significant accomplishment in the theory of approximation. This classical theorem states the following; {\it A sequence $\{L_n\}_{n\in\mathbb{N}}$ of linear positive operators on $C[0,1]$ satisfies $L_n(f)\rightarrow   f$  uniformly for every $f\in C[0,1]$ if and only if $L_n(e_j)\rightarrow   e_j$  uniformly for  $j=0,1,2,$ where $e_j(t)=t^j$ for $t\in[a,b]$ and $j=0,1,2$.} This theorem has various generalizations and analogues to different settings and many applications to various branches of science (see \cite{korovkin} for details). In $1968,$ O. Shisha and B. Mond  introduced the quantitative forms of Korovkin's theorem. These results include estimates for the convergence of the positive linear operators to the identity operator in terms of the convergence of the test functions and the quantity called modulus of continuity. The modulus of continuity is a measure of the continuity of the function, defined as follows;  
	\begin{definition}{\textbf{Modulus of continuity}}
		Let $f:[a,b]\rightarrow \mathbb{R}$ be a bounded function and let $\delta>0$. Then the modulus of continuity of $f$ with respect to the argument $\delta$ is defined as
		\[
		\omega(f,\delta)=\sup\{|f(x)-f(y)|:|x-y| \leq \delta, x,y\in [a,b]\}
		\]
	\end{definition}
  From the definition, it is evident that the modulus of continuity is always non-decreasing with respect to the argument $\delta$. Moreover any uniformly continuous function $f$ on $[a,b]$ satisfies $
  \lim\limits_{\delta\rightarrow 0^+}\omega(f,\delta)=0.$
Another important property is that for any $\lambda,\delta>0$, 
\[
\omega(f,\lambda\delta)\leq (1+[\lambda])\omega(f,\delta)
\]	
where $[\lambda]$	denotes the integer part of $\lambda$.
\begin{remark}
 Let $f\in C[a,b]$. $x,y\in[a,b]$ and $\delta>0$. For $|x-y|> \delta$, by the property of the modulus of continuity we have 
 \[
 |f(x)-f(y)|\leq \omega(f,|x-y|)\leq (1+|x-y|\delta^{-1})\omega(f,\delta)\leq (1+(x-y)^2\delta^{-2})\omega(f,\delta)\textrm{ for all }x,y\in [a,b].
 \]

\end{remark}

O. Shisha and B. Mond used  the above inequality  to prove the quantitative form of the Korovkin theorem in  \cite{shisha}
 	\begin{theorem}\cite{shisha}\label{shish}
 		Let $\{L_n\}\_{n\in\mathbb{N}}$ be a sequence of positive linear operators with the same domain $D$ which contains the restrictions of $1,t,t^2$
 		to $[a,b]\subseteq \mathbb{R}$. For $n=1,2,\ldots$, suppose that $L_n(1)$ is bounded. Let $f\in D$ be continuous in $[a,b]$, and $\|.\|$ the  sup-norm. For  $n=1,2,\ldots$, we have
 		\[
 		\|f-L_n(f)\|\leq\|f\|\|L_n(1)-1\|+ \|L_n(1)+1\|\omega(f,\mu_n)
 		\textrm{ where } 
 		\mu_n^2:=\|L_n((t-x)^2)(x)\|. 
 		\]
 		In particular, if $L_n(1)=1$, this inequality becomes 
 		$\|f-L_n(f)\|\leq 2\omega(f,\mu_n)$.
 	\end{theorem}

Recently, Dumitru Popa proved an operator version of the Korovkin theorem in  \cite{popa}. In this version, he considered sequence of positive linear maps between $C[a,b]$ and $C(T)$, where $T$ is a compact Hausdorff space, and an operator $A:C[a,b]\rightarrow C(T)$ with some special properties. We state the result below. $C(T)$ denotes the Banach space of all continuous functions $f:T\rightarrow \mathbb{R}$ with the supremum norm $\|f\|=\sup\limits_{t\in T}|f(t)|$.   
	\begin{theorem}\cite{popa}\label{popa}
		Let $T$ be a compact Hausdorff space, $V_n:C[a,b]\rightarrow C(T)$ a sequence of positive linear operators and $A:C[a,b]\rightarrow C(T)$ a linear operator such that $A(e_0)A(e_2)=[A(e_1)]^2$ and $A(e_0)(t)>0$ for every $ t\in T$. If $\lim\limits_{n\rightarrow\infty}V_n(e_j)=A(e_j)$ for $j=0,1,2$ all uniformly on $T$, then for every $f\in C[a,b]$, $\lim\limits_{n\rightarrow\infty}V_n(f)=A(f)$ uniformly on $T$.
		\end{theorem} 
As observed in \cite{popa},  this result applies to different examples, including those constructed from modified Kantorovich and Bernstein operators and their extensions.
	In this article, we prove a quantitative form of the above theorem, apply this result to various examples and obtain quantitative estimates for each. In some of these examples, we could also obtain the order of convergence. We also prove the trigonometric analogue of this result. Finally, we find an important application of our results to the preconditioners of linear systems with the Toeplitz structure. 
	
	\section{The Main Result}
	In this section, we prove an inequality involving modulus of continuity that gives a quantitative form of Theorem \ref{popa}. The modulus of continuity appears on the right side of this inequality. Therefore the rate of convergence is related to the modulus of continuity. This is precisely in the spirit of Theorem \ref{shish}.
\begin{theorem}\label{mains}
	Let $T$ be a compact Hausdorff space, $A:C[a,b]\rightarrow C(T)$ be a positive linear operator such that $A(e_0)(t)>0$  for every $ t\in T$. If $L_n:C[a,b]\rightarrow C(T)$ is a sequence of positive linear operators, then for every $f\in C[a,b]$, $n=1,2,\ldots$,

\begin{equation*}
\|L_n(f)-A(f)\| \leq 	m \{\|L_n(e_0)-A(e_0)\| \|A(f)\|+ (\|L_n(e_0)A(e_0)\|+ 1)\omega(f,\mu_n)\},
\end{equation*}
	where $ m=[(\inf\limits_{t\in T}A(e_0)(t))^{-1}],\, \mu_n^2=\|L_n(e_2)A(e_0)-2L_n(e_1)A(e_1)+ L_n(e_0)A(e_2)\|.$
	If, in addition, $L_n(e_0)=A(e_0)$, then the inequality becomes 
	\[
	\|L_n(f)-A(f)\| \leq m(\|A(e_0)^2\|+ 1)\omega(f,\mu_n).
	\]
	
\end{theorem}
\begin{remark}
	Note that
	\begin{multline*}
		L_n(e_2)A(e_0)-2L_n(e_1)A(e_1)+ L_n(e_0)A(e_2)= [L_n(e_2)-A(e_2)]A(e_0)-\\
		2[L_n(e_1)-A(e_1)]A(e_1)+[L_n(e_0)-A(e_0)]A(e_2) - 2\Delta,
	\end{multline*} 
	where $\Delta=A(e_0)A(e_2)-A(e_1)^2$. Therefore if $\Delta=0$ (this is part of assumptions in Theorem \ref{popa}), it is obvious that whenever $\lim\limits_{n\rightarrow\infty}L_n(e_j)=A(e_j)$ for $j=0,1,2$ uniformly on $T$, $\mu_n\rightarrow 0$ and we obtain that 
	$\lim\limits_{n\rightarrow \infty}L_n(f)=A(f)$ uniformly on $T$. Hence the operator version of Korovkin's theorem \cite{popa} is obtained as a consequence of our result.
\end{remark}
\begin{proof}[Proof of Theorem \ref{mains}]
	Let $x,y\in [a,b]$ and $\delta>0$. We recall the following inequality involving the modulus of continuity of $f$,
	$|f(x)-f(y)|\leq (1+(x-y)^2\delta^{-2})\omega(f,\delta) $
	That is for every $y\in [a,b]$ we have
$
	|f-f(y)|\leq (e_0 +(e_2-2ye_1+y^2e_0)\delta^{-2})\omega(f,\delta).
$
	
By the positivity of the operators $L_n, |L_n(f-f(y)e_0)|\leq L_n(|f-f(y)e_0|)$ which gives 
	\[
	|L_n(f)-f(y)L_n(e_0)|\leq (L_n(e_0)+ (L_n(e_2)-2yL_n(e_1)+ y^2L_n(e_0))\delta^{-2})\omega(f,\delta)
	\] 
	Hence for every $t\in T$ and $y\in[a,b]$,
	\begin{multline*}
		|L_n(f)(t)-f(y)L_n(e_0)(t)|\leq (L_n(e_0)(t)+ (L_n(e_2)(t)-2yL_n(e_1)(t)+\\ y^2L_n(e_0)(t))\delta^{-2})\omega(f,\delta)
	\end{multline*}
Now for each fixed $t\in T$,
	\begin{multline*}
		|L_n(f)(t)e_0-L_n(e_0)(t)f|\leq (L_n(e_0)(t)e_0+ (L_n(e_2)(t)e_0-2L_n(e_1)(t)e_1+\\ L_n(e_0)(t)e_2)\delta^{-2})\omega(f,\delta)
	\end{multline*}
Now we apply $A$. By the positive linearity of the operator $A$, for every $t\in T$,
	\begin{multline*}
		|L_n(f)(t)A(e_0)-L_n(e_0)(t)A(f)|\leq (L_n(e_0)(t)A(e_0)+ (L_n(e_2)(t)A(e_0)-\\
		2L_n(e_1)(t)A(e_1)+ L_n(e_0)(t)A(e_2))\delta^{-2})\omega(f,\delta)
	\end{multline*}
In particular for every $t\in T$,
	\begin{multline*}
		|L_n(f)(t)A(e_0)(t)-L_n(e_0)(t)A(f)(t)|\leq (L_n(e_0)(t)A(e_0)(t)+ (L_n(e_2)(t)A(e_0)(t)-\\
		2L_n(e_1)(t)A(e_1)(t)+ L_n(e_0)(t)A(e_2)(t))\delta^{-2})\omega(f,\delta)
	\end{multline*}
Taking supremum over $T$ we get,
	\begin{multline*}
		\|L_n(f)A(e_0)-L_n(e_0)A(f)\|\leq (\|L_n(e_0)A(e_0)\|+ (\|L_n(e_2)A(e_0)-
		2L_n(e_1)A(e_1)+ \\
		L_n(e_0)A(e_2)\|)\delta^{-2})\omega(f,\delta)
	\end{multline*}
	We choose $\delta^2=\mu_n^2=\|L_n(e_2)A(e_0)-2L_n(e_1)A(e_1)+ L_n(e_0)A(e_2)\|$. 
	We obtain
	\begin{equation*}
		\begin{split}
			\|L_n(f)A(e_0)-A(e_0)A(f)\|&\leq \|L_n(f)A(e_0)-L_n(e_0)A(f)\|+\|L_n(e_0)A(f)-A(e_0)A(f)\|\\
			&\leq (\|L_n(e_0)A(e_0)\|+ 1)\omega(f,\mu_n) + \|L_n(e_0)A(f)-A(e_0)A(f)\|
		\end{split}
	\end{equation*}
	Since $A(e_0)\in C(T)$, $A(e_0)(T)$ is a compact subset of $\mathbb{R}$. By the assumption $A(e_0)(t)>0$ for every $t \in T$, $\inf\limits_{t\in T} A(e_0)(t)=A(e_0)(t_0)$ for some $t_0 \in T$. Therefore $\inf\limits_{t\in T} A(e_0)(t)>0$, denote it by $\frac{1}{m}.$ Then we have  $\frac{1}{m}\|L_n(f)-A(f)\|\leq \|(L_n(f)-A(f))A(e_0)\|$ which completes the proof.
\end{proof}

\begin{remark}
	Suppose $\|L_n(e_2)A(e_0)-2L_n(e_1)A(e_1)+ L_n(e_0)A(e_2)\|=0$.  In the proof of Theorem \ref{mains}, for any $\delta>0$ we have, $\|L_n(f)A(e_0)-L_n(e_0)A(f)\|\leq \|L_n(e_0)A(e_0)\|\omega(f,\delta)$ and thus we obtain
\[
\|L_n(f)A(e_0)-A(e_0)A(f)\|\leq \|A(f)\|\|L_n(e_0)-A(e_0)\|+\|L_n(e_0)A(e_0)\|\omega(f,\delta)
\]
If $L_n=A$ for every $n\in \mathbb{N}$, then there is nothing to prove. So we avoid this trivial case. We let $\delta=\mu_n=\|L_n(e_0)-A(e_0)\|$. (If $\mu_n$ becomes $0$ then by the properties of the modulus of continuity we arrive at the trivial case where $L_n=A$ for every $n\in \mathbb{N}$). Hence we get a simplified estimate which depends only on the convergence of $L_n$ to $A$ on the single test function $e_0$ and the modulus of continuity of $f$ with argument $\|L_n(e_0)-A(e_0)\|$.
	\end{remark}

	\section{ Examples}
In this section, we discuss two  types of examples; Kantorovich and Bernstein type operators. These examples were discussed  in \cite{popa} as applications of Theorem \ref{popa}.  Here we obtain the convergence rate in each case. 
\subsection{Kantorovich-type operators}
\begin{example}\label{kant1}
	Let $T$ be a compact Hausdorff space, $(\Omega_n, m_n)_{n\in \mathbb{N}}$ be a sequence of probability measure spaces, $h_n:\Omega_n\rightarrow \mathbb{R}$ a sequence of $m_n$-integrable functions with the property that for every $n\in \mathbb{N}$, there exists $M_n\geq 0$ such that $0\leq h_n(\omega)\leq M_n$, $\forall \omega\in \Omega_n$. Let $a_n:T\rightarrow [0,1]$ be a sequence of continuous functions and $a:T\rightarrow [0,1]$, a continuous function. Define the operators $L_n, A:C[0,1]\rightarrow C(T)$ by
	\begin{equation*}
		L_n(f)(t)=\sum\limits_{k=0}^{n}\binom{n}{k}(a_n(t))^k(1-a_n(t))^{n-k}\int_{\Omega_n} f(\frac{k+h_n(\omega)}{n+M_n})dm_n(\omega);\hfill \,\,\,\,\,\,A(f)=f(\frac{a(t)+\beta}{\alpha+1}),
	\end{equation*}
where $\alpha=\lim\limits_{n\rightarrow\infty}\frac{M_n}{n}$, $\beta=\lim\limits_{n\rightarrow\infty}\frac{\int_{\Omega_n}h_n dm_n}{n}$, $\gamma=\lim\limits_{n\rightarrow\infty}\frac{\int_{\Omega_n}h_n^2 dm_n}{n^2}$ and assume that $\gamma=\beta^2$.

	From Theorem $2$ in \cite{popa}, we have the following;
$\lim\limits_{n\rightarrow\infty}L_n(f)=A(f)$ uniformly on $T$ for every $f\in C[0,1]$ if and only if $\lim\limits_{n\rightarrow\infty}a_n=a$ uniformly on $T$. We aim to obtain the convergence rate using Theorem \ref{mains}.

	Let $\frac{1}{1+\frac{M_n}{n}}=\frac{1}{\alpha+1}+ O(\alpha_n)$,  $\frac{\int_{\Omega_n}h_ndm_n}{n}=\beta+ O(\beta_n)$ and $\frac{\int_{\Omega_n}h_n^2dm_n}{n^2}=\beta^2+ O(\gamma_n)$ ($n\rightarrow \infty$).
	Using Theorem \ref{mains} and the facts that $L_n(e_0)=A(e_0)$ and $m=1$ here, we get
	\[\|L_n(f)-A(f)\|\leq 2\omega(f,\mu_n),
	\]
	where $\mu_n=O(\alpha_n+\beta_n+\gamma_n+\|a_n-a\|)$ ($n\rightarrow \infty$).

The following are some  variants of this example.
\begin{enumerate}
\item[ Case 1:](See Corollary $3$ of \cite{popa})

Let $c_n, d_n$ be sequence of real numbers satisfying $0\leq c_n\leq d_n$ for $n=1,2,\ldots$, $\alpha=\lim\limits_{n\rightarrow \infty}\frac{d_n}{n}$ and $\lim\limits_{n\rightarrow \infty}\frac{c_n}{n}=0.$  Setting $\Omega_n=[0,1]$, $h_n(x)=c_nx\leq d_n=M_n$,  $A(f)(t)=f(\frac{a(t)}{\alpha+1})$ for $t\in T$ and $L_n, a_n, a$ same as above, we get $\mu_n=O(\alpha_n+\|a_n-a\|)$ where $\frac{1}{1+\frac{d_n}{n}}=\frac{1}{1+\alpha}+O(\alpha_n)$.

\item[ Case 2:] (See Corollary $4$ of \cite{popa})

Consider $A(f)(t)=f(\frac{2a(t)+1}{4})$ for $t\in T, \Omega_n=[0,1]^n, h_n(x_1,\ldots,x_n)=x_1+\ldots+x_n\leq n=M_n$, $\alpha=1$.
We get $\alpha_n=\beta_n=0$ and $\gamma_n=\frac{1}{n},$ and hence $\mu_n^2=O(\frac{1}{n}+\|a_n-a\|)$ ($n\rightarrow \infty$).
\end{enumerate}

\end{example}

\begin{example}\label{kant2}
Define the operators $L_n, A:C[0,1]\rightarrow C(T)$ defined by 
	\[
	L_n(f)(t)=\sum\limits_{k=0}^{n}\binom{n}{k}(\frac{a_n(t)}{n})^k(1-\frac{a_n(t)}{2n})^{n-k}\int_{\Omega_n}f(\frac{k+h_n(\omega)}{n+M_n})dm_n(\omega);\,\,\,\,\,A(f)(t)=e^\frac{a(t)}{2}f(\frac{\beta}{\alpha+1}),
	\]	
where $\alpha$ $\beta$ and  $\gamma$ be as defined in Example \ref{kant1}. Then $A$ is a positive linear operator and $A(e_0)(t)=e^{\frac{a(t)}{2}}$ (this gives $m^{-1}>1$). Proposition $3$ of \cite{popa} states the following;  $\lim\limits_{n\rightarrow \infty}L_n(f)=A(f)$ uniformly on $T$  for $f\in C[0,1]$, if and only if $\lim\limits_{n\rightarrow \infty}a_n=a$ uniformly on $T$. Applying Theorem \ref{mains}, we get
	\[
	\|L_n(f)-A(f)\|\leq \|(1+\frac{a_n}{2n})^{n}-e^{\frac{a}{2}}\|\|A(f)\|+(\|(1+\frac{a_n}{2n})^{n}e^{\frac{a}{2}}\|+1)\omega(f,\mu_n)
	\]
	where $\mu_n^2=	O(\frac{1}{(1+\alpha)^2}\|e^{\frac{a}{2}}[\frac{a_n}{n^2}(1+\frac{a_n}{2n})^{n-1}+
	\frac{n-1}{n}(\frac{a_n}{n})^2(1+\frac{a_n}{2n})^{n-2}]\|+\alpha_n+\beta_n+\gamma_n) (n\rightarrow\infty)$ where $\alpha_n$, $\beta_n$ and $\gamma_n$ are as defined in Example \ref{kant1}.
	\end{example}

\subsection{Bernstein-type operators}
\begin{example}\label{bern1}
	Let $1\leq p<\infty,$ and $l_p=\{(x_n)_{n\in N}:\sum\limits_{n=1}^{\infty}|x_n|^p<\infty, x_n\in\mathbb{R}\}$ be the Banach space endowed with the $p$-norm $\|x\|_p=(\sum\limits_{n=1}^{\infty}|x_n|^p)^{\frac{1}{p}}$. For each $x=(x_n)_{n\in \mathbb{N}}\in l_p,$ denote $<x,e_n>:=x_n$, where $(e_n)_{n\in\mathbb{N}} $ denotes the standard Schauder basis in $l_p$. For a fixed $y=(y_n)_{n\in \mathbb{N}}\in l_p$ with $y_n\in [0,1]$, consider the compact subset $T_y$ of $l_p$ defined by $T_y=\{x=(x_n)_{n\in \mathbb{N}}\in l_p: 0\leq x_n\leq y_n\}$. Define $a_n, a:T_{y}\rightarrow [0,1]$ by $a_n(x)=<x,e_n>$  and $a(x)=0$. It is evident that $|a_n(x)|=O(y_n)$ as $n\rightarrow \infty$. Consider the sequence of positive linear operators $H_n, A:C[0,1]\rightarrow C(T_y)$ defined by 
	\[
	H_n(f)(x)=\sum_{k=0}^{n}\binom{n}{k}(<x,e_n>)^k(1-<x,e_n>)^{n-k}f(\frac{k}{n});\,\,\,\,\,\,A(f)(x)=f(0).
	\]
Corollary $5$, $(i)$ of \cite{popa} states that $\lim\limits_{n\rightarrow\infty} H_n(f)=A(f)$ uniformly in $T_y$. We apply Theorem \ref{mains}. and note that $m=1$. We have $\|H_n(f)-A(f)\|\leq 2\omega(f,\mu_n)$
	where $\mu_n=O(y_n)$ ($n\rightarrow \infty$) which is our desired quantitative estimate.\par
	Now if  $a_n, a$ are defined as 
	$a_n(x)=\sum\limits_{k=1}^{n}\frac{<x,e_k>}{2^k}$ and $a(x)=\sum\limits_{k=1}^{\infty}\frac{<x,e_k>}{2^k}$, then we have $\|a_n-a\|=O(2^{-n})$ ($n\rightarrow \infty$).  Now $G_n, A:C[0,1]\rightarrow C(T_y)$ be defined by
	\[
	G_n(f)(x)=\sum_{p=0}^{n}\binom{n}{p}(\sum\limits_{k=1}^{n}\frac{<x,e_k>}{2^k})^p(1-\sum\limits_{k=1}^{n}\frac{<x,e_k>}{2^k})^{n-p}f(\frac{p}{n});\,\,\,\,\,A(f)(x)=f(a(x)).
	\]
Corollary $5$, $(ii)$ of \cite{popa} states that for $f\in C[0,1]$, $\lim\limits_{n\rightarrow\infty} G_n(f)=A(f)$ uniformly in $T_y$. Here $m=1$ and we have the estimate $\|G_n(f)-A(f)\|\leq 2\omega(f,\mu_n)$
	where $\mu_n=O(2^{-n})$ ($n\rightarrow \infty$).
\end{example}

\begin{example}\label{bern2}
	Let $a_n,a:T\rightarrow [0,1]$ be continuous functions defined on a compact Hausdorff space $T$ and $L_n, A:C[0,1]\rightarrow C(T)$ be defined by 
	\[
	L_n(f)(t)=\sum\limits_{k=0}^{n}\binom{n}{k}(\frac{a_n(t)}{n})^k(1-\frac{a_n(t)}{2n})^{n-k}f(\frac{k}{n});\,\,\,\,\,\,A(f)=e^{\frac{a(t)}{2}}f(0).
	\]	
By Proposition $2$ of \cite{popa},  $\lim\limits_{n\rightarrow \infty}L_n(f)=A(f)$ uniformly on $T$ if and only if $\lim\limits_{n\rightarrow \infty}a_n=a$ uniformly on $T$. Applying Theorem \ref{mains}, and since $m^{-1}=\inf\limits_{t\in T}A(e_0)(t)=e^{\frac{a(t)}{2}}>1,\,a(t)\in[0,1]$ we get the following:
	\[
	\|L_n(f)-A(f)\|\leq \|A(f)\|\|(1+\frac{a_n}{2n})^n-e^{\frac{a}{2}}\| + (\|(1+\frac{a_n}{2n})^{n}\|+1)\omega(f,\mu_n)
	\]
	where $\mu_n^2=\|[\frac{a_n}{n^2}(1+\frac{a_n}{2n})^{n-1}+\frac {n-1}{n}.\frac{a_n^2}{n^2}(1+\frac{a_n}{2n})^{n-2}]e^{\frac{a}{2}}\|=O(\frac{1}{n^2})$.
\end{example}	
\begin{example}\label{bern3}
	Let $a:T\rightarrow [0,\infty)$ and $b:T\rightarrow (-\infty, 1]$ be continuous functions. Define $L_n, A:C[0,1]\rightarrow C(T)$ by 
	\[
	L_n(f)(t)=	\frac{1}{(a(t)+1)^n}\sum\limits_{k=0}^n\binom{n}{k}(a(t))^k(1-\frac{b(t)}{n})^{n-k}f(\frac{k}{n});\,\,\,\,A(f)=e^{-\frac{b(t)}{a(t)+1}}f(\frac{a(t)}{a(t)+1}).
	\]		
By proposition $4$ of \cite{popa},  $\lim\limits_{n\rightarrow \infty}L_n(f)=A(f)$ uniformly on $T$ for $f\in C[0,1]$.
Use Theorem \ref{mains}, and  $m^{-1}=\inf\limits_{t\in T}A(e_0)(t)>e^{-1}$ to get
	\[
	\|L_n(f)-A(f)\|\leq e\{\|(1-\frac{1}{n}\frac{b}{(a+1)})^n-e^{\frac{-b}{(a+1)}}\|+(\|(1-\frac{1}{n}\frac{b}{(a+1)})^ne^{\frac{-b}{(a+1)}}\|+1)\omega(f,\mu_n)\}
	\] 
	where $\mu_n^2=\frac{1}{n}\|e^\frac{-b}{a+1}(1-\frac{b}{n(a+1)})^{n-2}[(\frac{a}{a+1})^2(\frac{1}{n}(\frac{b}{a+1})^2-1)+a(1-\frac{b}{n(a+1)})]\|.$
\end{example}

\section{Trigonometric analogue and Applications}
The trigonometric analogue of the classical Korovkin theorem is well known and finds applications in various areas, including the convergence issues of the Fourier series. Here we prove the trigonometric analogue of Theorem \ref{popa} and apply this result to the preconditioning problems of large linear systems with Toeplitz structure. Recall that a Toeplitz linear system is related to the Fourier series expansion via the associated symbol-function. Hence the analogy with the applications of classical Korovkin theorems is evident. Here we directly prove the quantitative result first. The trigonometric analogue comes as a consequence.\par
$C_{2\pi}[-\pi,\pi]$ denotes the space of all continuous $2\pi$-periodic functions on the interval $[-\pi,\pi]$. Let $h_0(x)=1$, $h_1(x)=\cos x$ and $h_2(x)=\sin x$ for $x\in [-\pi,\pi]$. 
\begin{theorem}\label{maintri1}
	Let $T$ be a compact Hausdorff space, $A:C_{2\pi}[-\pi,\pi]\rightarrow C(T)$ be a positive linear operator such that $A(h_0)(t)>0$ for every $ t\in T$. If $L_n:C_{2\pi}[-\pi,\pi]\rightarrow C(T)$ is a sequence of positive linear operators, then for every $f\in C_{2\pi}[-\pi,\pi]$, $n=1,2,\ldots$,
	\begin{equation*}
		\|L_n(f)-A(f)\| \leq 	m \{\|L_n(h_0)-A(h_0)\| \|A(f)\|+ (\|L_n(h_0)A(h_0)\|+ 1)\omega(f,\mu_n)\},
	\end{equation*}
	where $ m=[(\inf\limits_{t\in T}A(h_0)(t))^{-1}],\, \mu_n^2=\frac{\pi^2}{2}\|L_n(h_0)A(h_0)-L_n(h_1)A(h_1)- L_n(h_2)A(h_2)\|.$
	If, in addition, $L_n(h_0)=A(h_0)$, then the inequality becomes 
	\[
	\|L_n(f)-A(f)\| \leq m(\|A(h_0)^2\|+ 1)\omega(f,\mu_n).
	\]
\end{theorem}

\begin{proof}
	Proof uses the same techniques of Theorem \ref{mains}, with some minor changes. For $f\in C_{2\pi}[-\pi,\pi], x,y\in [-\pi,\pi]$ and $\delta>0$, we have
	\begin{equation*}
		\begin{split}
	|f(x)-f(y)|&\leq (1+(x-y)^2\delta^{-2})\omega(f,\delta)\\
	&\leq (1+\pi^2\sin^2\frac{x-y}{2}\delta^{-2})\omega(f,\delta)\\
	&=(1+\frac{\pi^2}{2}(1-\cos(x-y))\delta^{-2})\omega(f,\delta)\\
	&=(1+\frac{\pi^2}{2}(1-\cos x\cos y-\sin x\sin y)\delta^{-2})\omega(f,\delta)\\
	&=(h_0(x)h_0(y)+\frac{\pi^2}{2}(h_0(x)h_0(y)-h_1(x)h_1(y)-h_2(x)h_2(y))\delta^{-2})\omega(f,\delta))
\end{split}
\end{equation*} 
	Hence for any $y\in [-\pi,\pi]$ we have\\
	$|f-f(y)|\leq (h_0.h_0(y)+\frac{\pi^2}{2}(h_0.h_0(y)-h_1.h_1(y)-h_2.h_2(y))\delta^{-2})\omega(f,\delta).$
Now proceeding in a similar fashion as in Theorem \ref{mains}, we obtain
\begin{multline*}
	\|L_n(f)A(h_0)-L_n(h_0)A(f)\|\leq (\|L_n(h_0)A(h_0)\|+ \\
	\frac{\pi^2}{2}(\|L_n(h_0)A(h_0)-L_n(h_1)A(h_1)-L_n(h_2)A(h_2)\|)\delta^{-2})\omega(f,\delta).
\end{multline*} 
Choose $\delta^2=\mu_n^2=\frac{\pi^2}{2}\|L_n(h_0)A(h_0)-L_n(h_1)A(h_1)-L_n(h_2)A(h_2)\|$ to complete the proof.
    
\end{proof}
\begin{remark}
	$L_n(h_0)A(h_0)-L_n(h_1)A(h_1)-L_n(h_2)A(h_2)=(L_n(h_0)-A(h_0))A(h_0)-(L_n(h_1)-A(h_1))-(L_n(h_2)-A(h_2))A(h_2)+ [A(h_0)^2)-(A(h_1)^2+A(h_2)^2)]$. 
	Note that if $A(h_0)^2=A(h_1)^2+A(h_2)^2$, then $\lim\limits_{n\rightarrow\infty}L_n(h_i)=A(h_i)$ uniformly on $[-\pi,\pi]$ for $i=0,1,2$ gives $\lim\limits_{n\rightarrow\infty}L_n(f)=A(f)$ uniformly on $[-\pi,\pi]$ for every $f\in C_{2\pi}[-\pi,\pi]$.
\end{remark}
As a consequence of the above remark, we obtain the following trigonometric analogue of the operator version of the Korovkin theorem.
\begin{theorem}\label{maintri2}
	Let $T$ be a compact Hausdorff space, $V_n:C_{2\pi}[-\pi,\pi]\rightarrow C(T)$ a sequence of positive linear operators and $A:C_{2\pi}[-\pi,\pi]\rightarrow C(T)$ a linear operator such $A(h_1)^2+A(h_2)^2=[A(h_0)]^2$ and $A(h_0)(t)>0$ for every $ t\in T$. If $\lim\limits_{n\rightarrow\infty}V_n(h_i)=A(h_i)$ for $i=0,1,2$ all uniformly on $T$, then for every $f\in C_{2\pi}[-\pi,\pi]$, $\lim\limits_{n\rightarrow\infty}V_n(f)=A(f)$ uniformly on $T$.	
\end{theorem}
\begin{remark}
It is worthwhile to notice that the conditions $A(h_1)^2+A(h_2)^2=A(h_0)^2$ and $A(h_0)(t)>0$ are quite natural. For, in the classical trigonometric analogue, the test functions are $1, \sin x$ and $ \cos x$, and they satisfy the conditions, with $A$ being the identity operator. 
\end{remark}
\subsection{Application to Preconditioners of Toeplitz systems}

Now we briefly describe the preconditioners of Toeplitz systems generated by continuous real-valued functions and discuss a positive linear operator associated with them. We apply our theorem and obtain quantitative estimates for such operators.\par
Let $M_n(\mathbb{C})$ be the collection of all $n\times n$ complex matrices  equipped with the Frobenius inner-product defined by $<A,B>_F=trace(B^*A)$ ($*$ denotes conjugate transpose). $M_n(\mathbb{C})$ is a  Hilbert space with respect to the Frobenius norm $\|A\|_F=(\sum\limits_{k=0}^{n-1}|A_{i,j}|^2)^\frac{1}{2}$ induced by the inner-product.
Let $U_n\in M_n(\mathbb{C})$ be a unitary matrix. We define
\[
M(U_n)=\{A\in M_n(\mathbb{C}):A=U_nDU_n^*\ where \ D\ is\ a\ complex\ diagonal\ matrix.\}
\]
Then $M(U_n)$ is a closed sub-algebra of $M_n(\mathbb{C})$. For each $A\in M_n(\mathbb{C})$, 
there exists a unique best approximation $P_{U_n}(A)\in M(U_n)$ to $A$ such that $\|P_{U_n}(A)-A\|_F= \min\{\|A-B\|_F:B\in M(U_n)\}$. We call $P_{U_n}(A)$ "preconditioner" of $A$.\par
Preconditioners have applications in solving complicated linear systems of the form $Ax=b$ by replacing them with simpler ones and obtaining approximate solutions. This is a widely used technique in numerical linear algebra to address the stability issues during the iteration processes and to convert an ill-posed system into a well-posed one. It is a simple exercise ( see \cite{stefano} for example) to show that $P_{U_n}(A)=U_n\delta(U_n^*AU_n)U_n^*$ where $\delta(A)$ denotes the diagonal matrix obtained by setting the off-diagonal elements of $A$ to $0$ and diagonal entries same as that of $A$. For a detailed description of preconditioners and for the Korovkin-type theorems in this direction, we refer to \cite{stefano, studia13, laa18, survey21}. We consider the case when $A=A_n$ is a sequence of Toeplitz matrices and $U_n$ is a particular unitary matrix sequence.\par
Let $f$ be a continuous $2\pi$-periodic real valued function defined on a real interval $I$.
The $n\times n$ Toeplitz matrix $A_n(f)$ generated by $f$ is defined by
\[
(A_n(f))_{i,j}=\hat{f}(i-j),
\]
where $\hat{f}(m)$ denotes the $m^{th}$ Fourier coefficient of $f$, defined by $
\hat{f}(m)=\frac{1}{2\pi}\int_{-\pi}^{\pi}f(x)e^{-imx}\ dx.$

In the rest of this article, we consider only Toeplitz matrix sequences. Circulant matrices are a special type of Toeplitz matrices, where the whole matrix is obtained by simply shifting the first row. These matrices are characterized by their  (unitary) diagonalizability with respect to the  Fourier basis. The corresponding sequence of unitary matrices are given by $U_n=(\frac{1}{\sqrt{n}}e^{ij\frac{2\pi k}{n}})$, $j,k=0,\ldots,n-1$.  The usage of circulant preconditioners to the Toeplitz linear system is a folklore idea in numerical linear algebra. This idea is well established theoretically and supported by efficient algorithms due to Tony Chan, R. H. Chan, G. Strang, etc. (see \cite{Tony, ChanRH, strang} and references therein). We consider the generalized Vandermonde matrices as  $U_n$, covering the circulants as a special case and including more efficient preconditioners such as Hartley, Wavelets, etc.

Let for each $n\in\mathbb{N}$, $\{v_{n,j}\}_{j=0}^{n-1}$ be a finite sequence trigonometric functions on a real interval $I$. Let $G_n=\{x_k^{(n)}:k=0,1,\ldots,n-1\}$ be a sequence of $n$ grid points on $I$. Consider the generalized Vandermonde matrix $V_n=(v_{n,j}({x_k^{(n)})})_{k,j=0}^{n}$ and assume that $U_n:=V_n^*$ is a unitary matrix. Then $M(U_n)$ forms a trigonometric matrix algebra. Define a sequence of linear operators $L_n[U_n]$ on $C_{2\pi}(I)$ by 
\[
L_n[U_n](f)(x)=v(x)A_n(f)v^*(x),\ x\in I.
\] 
\begin{remark}
The importance of the linear maps  $L_n[U_n]$ lies in the fact that the eigenvalues of the preconditioners $P_{U_n}(A_n(f)))$ are obtained as evaluation of  $L_n[U_n](f)$ at some grid points in $I.$  Since the spectral information of the Toeplitz matrices is stored in the symbol function $f,$ and the spectral information of the preconditioners is stored in $L_n[U_n](f)$, it is important to know if $L_n[U_n](f)\rightarrow f$, to understand the asymptotic spectral behaviour of $A_n(f)$.
\end{remark}
It is well known from \cite{stefano} that $L_n[U_n]$ is a positive linear operator. We apply Theorem \ref{maintri1}, \ref{maintri2} to $L_n[U_n]$ and obtain convergence rates.
\begin{example}
Let $L_n[U_n], A:C_{2\pi}[-\pi,\pi]\rightarrow C_{2\pi}[-\pi,\pi]$ be defined by 
$
L_n[U_n](f)(x)=v(x)A_n(f)v^*(x),$ and  $A(f)(x)=v(x)v^*(x)f(x)$, $x\in [-\pi,\pi]$. Assume that $v(x)\neq 0$ for every $x\in [-\pi,\pi]$. Then $A(h_0)(x)=v(x)v^*(x)>0$ (Note that $v(x)v^*(x)=0 \implies v(x)=0$). Note that here $A(h_0)^2=(v(x)v^*(x))^2.1=(v(x)v^*(x))^2h_1^2+ (v(x)v^*(x))^2h_2^2=A(h_1)^2+A(h_2)^2$. Thus Theorem \ref{maintri2} is applicable here.

Since $A_n(h_0)=I_n$, the $n\times n$ identity matrix, the estimate in Theorem \ref{maintri1}, becomes 
\[
\|L_n(f)-A(f)\|\leq 2m\omega(f,\mu_n),
\]
where $\mu_n^2=\frac{\pi^2}{2}\|L_n(h_0)A(h_0)-L_n(h_1)A(h_1)-L_n(h_2)A(h_2)\|$. 
\end{example}

\begin{example}
	Let $a:[-\pi,\pi]\rightarrow [0,1]$ be a continuous $2\pi$-periodic function. Consider the sequence of positive linear operators  
	\[
	L_n[U_n](f)(x)=v(x)A_n(e^{a(x)}f)v^*(x)=e^{a(x)}v(x)A_n(f)v^*(x)
	\]
	where $U_n=(\frac{1}{\sqrt{n}}e^{ijx_k^{(n)}})$, $k,j=0,\ldots,n-1$ and the grid points $G_n=\{x_k^{(n)}:\frac{2\pi k}{n}:k=0,\ldots,n-1\}\subset [-\pi,\pi]$. 
	Hence the vector $v(x)=\frac{1}{\sqrt{n}}(1, e^{ix}, e^{2ix},..., e^{(n-1)ix})$. Let $A(f)(x)=e^{a(x)}f(x)$, $x\in[-\pi,\pi]$. We have $A(h_0)(x)=e^{a(x)}>0$ and $m\leq1$. $A(h_1)(x)=e^{a(x)}\cos x$ and $A(h_2)(x)=e^{a(x)}\sin x$ and hence $A(h_0)^2=A(h_1)^2+A(h_2)^2$.\\
		Also, it is easy to see that
	$A_n(h_0)=I$,\\
	$A_n(h_1)=
	\begin{pmatrix}
		0&\dfrac{1}{2}&0&\ldots&0&\\
		\dfrac{1}{2}&0&\dfrac{1}{2}&\ldots&0&\\
		\ldots&\dfrac{1}{2}&\ddots&\ddots&\ldots&\\
		\ldots&\ldots&\ddots&\ddots&\dfrac{1}{2}&\\
		0&\ldots&0&\dfrac{1}{2}&0
	\end{pmatrix}$,
	$A_n(h_2)=
	\begin{pmatrix}
		0&\dfrac{i}{2}&0&\ldots&0&\\
		\dfrac{-i}{2}&0&\dfrac{i}{2}&\ldots&0&\\
		\ldots&\dfrac{-i}{2}&\ddots&\ddots&\ldots&\\
		\ldots&\ldots&\ddots&\ddots&\dfrac{i}{2}&\\
		0&\ldots&0&\dfrac{-i}{2}&0
	\end{pmatrix}$.\\
	$L_n[U_n](h_0)(x)=e^{a(x)}$, $L_n[U_n](h_1)(x)= \frac{n-1}{n}e^{a(x)}\cos x$, $L_n[U_n](h_2)=v(x)A_n(h_2)v^*(x)=\frac{n-1}{n}e^{a(x)}\sin x$.
$|L_n(h_0)(x)A(h_0)(x)-L_n(h_1)(x)A(h_1)(x)-L_n(h_2)(x)A(h_2)(x)|\\
=e^{2a(x)}|1-\frac{n-1}{n}\cos^2 x-\frac{n-1}{n}\sin^2 x|=\frac{e^{2a(x)}}{n}$.\\
	 Hence we get the estimate
	\[
	 \|L_n(f)-A(f)\|\leq (\|e^{2a}\|+1)\omega(f,\mu_n),
	\]
	where $\mu_n^2=\frac{\pi^2}{2}\frac{\|e^{2a}\|}{n}$, for $n=1,2,\ldots $\par
	If we consider the special case where $a(x)=0$, we get  $\mu_n^2=\frac{\pi^2}{2}\frac{1}{n}$
    and	\[
    \|L_n(f)-A(f)\|\leq 2\omega(f,\frac{\pi}{\sqrt{2n}}).
    \]
\end{example}
\begin{example}\label{fina1}
	Let $a_n, a:[-\pi,\pi]\rightarrow [0,1]$ be continuous $2\pi$-periodic functions. Consider the sequence of positive linear operators  
	\[
	L_n[U_n](f)(x)=e^{a_n(x)}v(x)A_n(f)v^*(x),
	\] 
	and take  $U_n=(\frac{1}{\sqrt{n}}e^{ijx_k^{(n)}})$, $k,j=0,\ldots,n-1$ with the grid points $x_k^{(n)}=\frac{2\pi k}{n}:k=0,\ldots,n-1.$\par
	Take $A(f)(x)=e^{a(x)}f(x)$, $x\in[-\pi,\pi]$. We have $A(h_0)(x)=e^{a(x)}>0$, $m<1$, $A(h_1)(x)=e^{a(x)}\cos x$ and $A(h_2)(x)=e^{a(x)}\sin x$ and $A(h_0)^2=A(h_1)^2+A(h_2)^2$.

		$L_n[U_n](h_0)(x)=e^{a_n(x)}$, $L_n[U_n](h_1)(x)= \frac{n-1}{n}e^{a_n(x)}\cos x$,  $L_n[U_n](h_2)=\frac{n-1}{n}e^{a_n(x)}\sin x$.\\
	$|L_n(h_0)(x)A(h_0)(x)-L_n(h_1)(x)A(h_1)(x)-L_n(h_2)(x)A(h_2)(x)|
	=\frac{e^{a(x)+a_n(x)}}{n}$\\
	We get the estimate
	\[
	\|L_n(f)-A(f)\|\leq \|e^{a_n}-e^{a}\|\|e^{a}f\|+(\|e^{a_n+a}\|+1)\omega(f,\mu_n)
	\]
	where $\mu_n^2=\frac{\pi^2}{2}\frac{e^{a+a_n}}{n}$. If we also have $\lim\limits_{n\rightarrow\infty}a_n=a$ uniformly on $[-\pi,\pi]$, then by Theorem \ref{maintri2}, we obtain
	$\lim\limits_{n\rightarrow\infty}L_n(f)=A(f)$ uniformly on $[-\pi,\pi]$ for every $f\in C_{2\pi}[-\pi,\pi]$.
\end{example}
\begin{example}
	Under the same conditions in Example \ref{fina1}, if we consider the operators $
	L_n[U_n](f)(t)=v(a_n(t))A_n(f)v^*(a_n(t))$ and $A(f)(t)=f(a(t))$, we get  $A(h_0)(t)=1$, $m=1$, $A(h_1)(t)=\cos a(t)$ and $A(h_2)(t)=\sin a(t)$ and $A(h_0)^2=A(h_1)^2+A(h_2)^2$.
		Also, $L_n[U_n](h_0)(t)=1$,
	$L_n[U_n](h_1)(t)=\frac{n-1}{n}\cos a_n(t)$,  $L_n[U_n](h_2)(t)=\frac{n-1}{n}\sin a_n(t)$.\\
	$|L_n(h_0)(t)A(h_0)(t)-L_n(h_1)(t)A(h_1)(t)-L_n(h_2)(t)A(h_2)(t)|
	=|1-\frac{n-1}{n}\cos (a_n(t)-a(t))|$\\
		We get the estimate $
	\|L_n(f)-A(f)\|\leq 2\omega(f,\mu_n),$
	where $\mu_n^2=\frac{\pi^2}{2}\|1-\frac{n-1}{n}\cos (a_n-a)\|.$
\end{example}
\section*{Concluding Remarks}
We restricted our study in the example section to the circulants case only as we get the order of convergence in this case. Many other preconditioners like Hartley, Tau, wavelets etc. can also be considered similarly. We can also consider extending these examples to the classical Jackson, Fejer type operators, etc. Whether the operator version has a non-commutative analogue is an interesting and challenging question to address in the future.
\section*{Acknowledgement}
V. B. Kiran Kumar is supported by the KSYSA- Research Grant by KSCSTE, Kerala. P.C.Vinaya is thankful to the University Grants Commission (UGC) for financial support.

	\nocite{*}
\bibliographystyle{elsarticle-num} 
\bibliography{Operatorversion}


	
	
	
\end{document}